\documentclass[12pt]{amsart}
\usepackage{amsmath,amssymb}
\newtheorem{theorem}{Theorem}
\newtheorem{proposition}[theorem]{Proposition}
\newtheorem{lemma}[theorem]{Lemma}

\def\C{\Bbb C}
\def\D{\Bbb D}
\def\O{\mathcal O}
\def\eps{\varepsilon}

\title[Comparison of invariant functions]{Comparison of invariant functions
on strongly pseudoconvex domains}

\author{Nikolai Nikolov}

\address{Institute of Mathematics and Informatics\\ Bulgarian Academy
of Sciences\\ Acad. G. Bonchev 8, 1113 Sofia, Bulgaria\newline
\indent Faculty of Information Sciences\\
State University of Library Studies and Information Technologies\\
Shipchenski prohod 69A, 1574 Sofia, Bulgaria}\email{nik@math.bas.bg}



\begin{document}

\begin{thanks}{The author would like to thank the referee for her/his remarks on
the paper.}
\end{thanks}

\begin{abstract}{It is shown that the Carath\'eodory distance and the Lempert
function are almost the same on any strongly pseudoconvex domain
in $\C^n;$ in addition, if the boundary is $C^{2+\eps}$-smooth, then
$\sqrt{n+1}$ times one of them almost coincides with
the Bergman distance.}
\end{abstract}

\maketitle

\section{Introduction and results}

Let $D$ be a bounded domain in $\C^n.$ Recall that the Bergman
distance $b_D$ is the integrated form of the Bergman metric
$\beta_D,$ i.e.
$$b_D(z,w)=\inf_\gamma\int_0^1\beta_D(\gamma(t);\gamma'(t))dt,$$
where the infimum is taken over all smooth curves $\gamma:[0,1]\to
D$ with $\gamma(0)=z$ and $\gamma(1)=w.$ Note that
$$\beta_D(z;X)=\frac{M_D(z;X)}{K_D(z)}, \quad z\in D,\  X\in\C^n,$$
where
$$M_D(z;X)=\sup\{|f'(z)X|:f\in L_h^2(D),\;||f||_D=1,\;f(z)=0\}$$ and
$$K_D(z)=\sup\{|f(z)|:f\in L_h^2(D),\;||f||_D=1\}$$ is the square root of the
Bergman kernel on the diagonal.

The Kobayashi distance $k_D$ is the largest distance not exceeding
the Lempert function
$$l_D(z,w)=\inf\{\tanh^{-1}|\alpha|:\exists\varphi\in\O(\D,D)
\hbox{ with }\varphi(0)=z,\varphi(\alpha)=w\},$$ where $\D$ is the
unit disc. Note that $k_D$ is  the integrated form of the
Kobayashi metric
$$\kappa_D(z;X)=\inf\{|\alpha|:\exists\varphi\in\O(\D,D) \hbox{
with }\varphi(0)=z,\alpha\varphi'(0)=X\}.$$

Finally, the Carath\'eodory distance is given by
$$c_D(z,w)=\sup\{\tanh^{-1}|f(w)|:f\in\O(D,\D)
\hbox{ with }f(z)=0\}.$$ Note that $c_D\le k_D\le l_D$ and $c_D\le b_D.$

Let now $\partial D$ be $C^2$-smooth. Recall that $D$ is called
strongly pseudoconvex if for any $p\in\partial D$ the Levi form of
the signed distance to $\partial D$ is positive definite on the
tangent space $T_p^\C\partial D.$

Z. M. Balogh and M. Bonk \cite{BB1,BB2} introduce a positive
function $g_D$ on $D\times D$ which depends on
$d_D=\mbox{dist}(\cdot,\partial D)$ and the Carnot-Carath\'eodory
metric on $\partial D.$ Estimating $\kappa_D,$ they prove that the
function $g_D-k_D$ is bounded on $D\times D$ (see \cite[Corollary
1.3]{BB2}). If $D\subset\C^n$ is $C^3$-smooth, then the same
result for $c_D$ and $b_D/\sqrt{n+1}$ is announced in \cite{BB1}
without proof  (see also \cite{BB2}); this means that the differences
$k_D-c_D$ and $k_D-b_D/\sqrt{n+1}$ are bounded on $D\times D.$

Recall another result in this direction due to S. Venturini
\cite[Theorem 1]{Ven}:

If $D$ is a strongly pseudoconvex domain in $\Bbb C^n,$ then
\begin{equation}\label{v}
\lim_{\substack{z\to\partial D\\z\neq
w}}\frac{c_D(z,w)}{k_D(z,w)}=1\quad\mbox{uniformly in }w\in D.
\end{equation}

Our first observation can be considered as a version of Lempert's
theorem (see \cite[Theorem 1]{Lem}), saying that $l_D=c_D$ on any
convex domain $D.$ In fact, this observation tells us that the
Balogh-Bonk results hold even for $l_D$ and $c_D,$ i.e. these two
functions are almost the same on any strongly pseudoconvex domain in
$\C^n.$

\begin{proposition}\label{1} If $D$ is a strongly pseudoconvex
domain in $\Bbb C^n,$ then the difference $l_D-c_D$ is bounded on
$D\times D.$
\end{proposition}

\noindent{\it Proof.} It is enough to show that
$$\limsup_{z\to p,w\to q}(l_D(z,w)-c_D(z,w))<\infty,$$
where $p\in\partial D$ and $q\in\overline{D}.$

If $p\neq q,$ this a consequence of the estimate
\begin{equation}\label{n}
2l_D(z,w)\le-\log d_D(z)-\log d_D(w)+c,\ z,w\in D
\end{equation}
(see \cite[Theorem 1]{NPT}) and the Vormoor-Fadlalla-Abate estimate
\begin{equation}\label{m}2c_D(z,w)\ge-\log d_D(z)-\log d_D(w)-c,\ z\hbox{ near p,}\
w\hbox{ near q}
\end{equation}(cf. \cite[Theorem 10.2.1 and Remark
10.2.2 (b)]{JP}), where $c>0$ is a constant.

When $p=q,$ we shall use that, according to the Fornaess embedding
theorem (see \cite[Proposition 1]{For}), there exist a holomorphic
map $\Phi$ from a neighborhood of $\overline{D}$ to $\C^n$ and a
strictly convex domain $G\supset\Phi(D)$ such that, near $p,$ $\Phi$
is biholomorphic and $\partial G=\partial\Phi(D).$

Now we may choose a neighborhood $U$ of $p$ such that $\Phi$ is biholomorphic
on $U$ and $G':=\Phi(D\cap U)$ is a strictly convex domain.
Set $z'=\Phi(z)$ and $w'=\Phi(w).$ Then, by Lempert's theorem and
the Balogh-Bonk result,
$$l_D(z,w)\le l_{D\cap U}(z,w)=l_{G'}(z',w')
\le g_{G'}(z',w')+c=g_G(z',w')+c$$
$$\le l_G(z',w')+2c=c_G(z',w')+2c\le c_D(z,w)+2c,\ z,w\hbox{ near }p.\qed$$

A similar result to Proposition \ref{1} holds for $b_D/\sqrt{n+1}$
instead of $c_D,$ i.e. $l_D$ and $b_D/\sqrt{n+1}$ almost coincide on
any strongly pseudoconvex domain in $\C^n.$

\begin{proposition}\label{2} If $D$ is a $C^{2+\eps}$-smooth strongly
pseudoconvex domain in $\C^n,$ then the difference
$l_D-b_D/\sqrt{n+1}$ is bounded on $D\times D.$
\end{proposition}

\noindent{\it Proof.} We have to show the boundedness of the
difference $b_D/\sqrt{n+1}-g_D.$ By \cite[Theorem 1.1]{BB2}, this is
a consequence of the following estimate for $\beta_D=M_D/K_D:$
\begin{equation}\label{s}e^{-Cd_D^s(z)}A_D(z;X)\le\beta_D(z;X)
\le e^{Cd_D^s(z)}A_D(z;X),
\end{equation}
where
$$A_D(z;X)=\frac{||X_n||}{4d_D^2(z)}+\frac{L_{q_z}(X_t)}{d_D(z)},$$
$s=\min(\eps/2,1/2),$ $C>0$ is a constant, $z$ is near $\partial D,$
$q_z\in\partial D$ is the closest point to $z,$ $L_{q_z}$ is the
Levi form of $-d_D$ at $q_z,$ and $X_n$ and $X_t$ are the complex
normal and tangential components of $X$ at $q_z,$ respectively.

The estimate \eqref{s} follows by the next estimates for $K_D$ and
$M_D:$
$$e^{-Cd_D^s(z)}B_D(z)\le
K^2_D(z)e^{Cd_D^s(z)}B_D(z),$$
$$e^{-Cd_D^s(z)}A_D(z;X)B_D(z)\le
M_D(z;X)\le e^{Cd_D^s(z)}A_D(z;X)B_D(z),$$ where
$$B_D(z)=\frac{n!}{4\pi^n}.\frac{\Pi_{q_z}}{d^{n+1}_D(z)},$$ and
$\Pi_{q_z}$ is the product of the eigenvalues of $L_{q_z}.$

Theorem B in \cite{Ma} gives the estimate \eqref{s} for $\kappa_D$
instead of $\beta_D$ with $s=1/2$ if $\partial D$ is $C^3$-smooth
(see also \cite[Proposition 1.2]{BB2} for a weaker estimate in the
$C^2$-smooth case). Since $K_D$ and $M_D$ are monotone under
inclusion of domains in $\C^n,$ the above estimates for these
invariants can be obtained in the same way as \cite[Theorems A and
B]{Ma} by two changes. The first one is the replacing of the
localization \cite[(2.6)]{Ma} for $\kappa_D$ (see also
\cite[(4.10)]{BB2}) by Lemma 3 below (in fact, there the weaker
inequalities $\ge 1+cd_D^s(z)$ are enough). The
$C^{2+\eps}$-smoothness instead of $C^3$- reflects to the second
change: the exponent $3/2$ in the inequality $|M(w)|\le
C|w_1|^{3/2}$ on the bottom of \cite[p.~334]{Ma} is replaced by
$1+s$ and then the inequality $\bigl| ||\Psi(w)||^2-1\bigr|\le
C\sqrt u$ on the top of \cite[p.~335]{Ma} becomes $\bigl|
||\Psi(w)||-1\bigr|\le Cu^s.$

\begin{lemma}\label{3} Let $q$ be a strongly pseudoconvex boundary
point of a pseudoconvex domain $D$ (not necessary bounded) and let
$U$ be a neighborhood of $q.$ There exist a neighborhood $V\subset
U$ of $q$ and a constant $c>0$ such that if $z\in D\cap V$ and
$X\in(\C^n)_\ast,$ then
$$\frac{K_D(z)}{K_{D\cap U}(z)}\ge 1+cd_D(z)\log d_D(z),$$
$$\frac{M_D(z;X)}{M_{D\cap U}(z;X)}\ge 1+cd_D(z)\log d_D(z).$$
\end{lemma}

\noindent{\it Proof.} Since $D$ is locally strongly convexifiable,
shrinking $U$ (if necessary), we may assume that for any $\tilde
q\in\partial D\cap U$ there exists a peak function $p_{\tilde q}$
for $D\cap U$ at $\tilde q$ such that $|1-p_{\tilde q}|\le
c_1||\cdot-\tilde q||,$ where the constant $c_1>0$ does not depend
on $\tilde q.$ Take neighborhoods $V\subset V_1\Subset V_2\Subset
U$ of $p$ with the following property: if $z\in D\cap V,$
$q_z\in\partial D$ and $||q_z-z||=d_D(z),$ then $q_z\in V_1.$ Note
that $1>c_2:=\sup_{z\in D\cap V}\sup_{D\cap U\setminus
V_2}|p_{q_z}|.$ Using H\"ormander's $L^2$-estimates for the
$\overline\partial$ operator (see \cite [Theorem
2.2.1$^\prime$]{Hor}) provides a constant $c_3>0$ such that for
any $k\in\Bbb N,$ $z\in D\cap V$ and $X\in(\C^n)_\ast$ one has
$$\frac{K_D(z)}{K_{D\cap
U}(z)}\ge\frac{|p_{q_z}(z)|^k}{1+c_3c_2^k},\quad
\frac{M_D(z;X)}{M_{D\cap U}(z;X)}
\ge\frac{|p_{q_z}(z)|^k}{1+c_3c_2^k}$$ (cf. the proof of
\cite[Theorem 2]{Nik} which is a variation of the proof of
\cite[Lemma 3.5.2]{Hor}).

We may assume that $d_D<1$ on $D\cap V.$ Since $|p_{q_z}(z)|^k\ge
1-c_1kd_D(z),$ choosing $k$ to be the integer part of $\log
d_D(z)/\log c_2,$ we find a constant $c>0$ with the desired
property.\qed
\smallskip

Finally, we point out the following

\begin{proposition}\label{4} Let $D$ be a strongly pseudoconvex
domain in $\C^n.$ Then
\begin{equation}\label{g}
\lim_{\substack{z\to\partial
D\\z\neq w}}\frac{c_D(z,w)}{l_D(z,w)}=1\quad\mbox{uniformly in }w\in
D.
\end{equation}
and
\begin{equation}\label{b} \lim_{\substack{z\to\partial D\\z\neq
w}}\frac{b_D(z,w)}{k_D(z,w)}=\sqrt{n+1}\quad\mbox{uniformly in }w\in
D.
\end{equation}
\end{proposition}

Note that, since $c_D\le k_D\le l_D,$ \eqref{g} is a generalization
of \eqref{v}. On the other hand, the proof of \eqref{g} uses \eqref{v}.
\smallskip

\noindent{\it Proof of \eqref{g}.} It is enough to show that
$$\liminf_{\substack{z\to p\\w\to q\\z\neq
w}}\frac{c_D(z,w)}{l_D(z,w)}\ge 1,$$ where $p\in\partial D$ and
$q\in\overline{D}.$

If $p\neq q,$ this is a consequence of the inequalities \eqref{n}
and \eqref{m}.

Let $p=q.$ By \cite[Proposition 3]{Ven}, we have that
\begin{equation}\label{t}
\liminf_{\substack{z\to p\\w\to p\\z\neq w}}\frac{k_{D\cap
U}(z,w)}{k_D(z,w)}=1
\end{equation} for any neighborhood $U$ of $p$
such that $D\cap U$ is connected. Choosing $U$ such that $D\cap U$
to be convexifiable, \eqref{v}, \eqref{t} and Lempert's theorem
imply that
$$\liminf_{\substack{z\to p\\w\to q\\z\neq
w}}\frac{c_D(z,w)}{l_D(z,w)}=\liminf_{\substack{z\to p\\w\to
q\\z\neq w}}\frac{k_D(z,w)}{l_D(z,w)}\ge \liminf_{\substack{z\to
p\\w\to q\\z\neq w}}\frac{k_{D\cap U}(z,w)}{l_{D\cap U}(z,w)}=1.$$
\smallskip

\noindent{\it Proof of \eqref{b}.} Due to K. Diederich and D. Ma
(see e.g. \cite[Theorems 10.4.2 and 10.4.6]{JP}), we have that
$$\lim_{z\to\partial D}\frac{\beta_D(z;X)}{\kappa_D(z;X)}=\sqrt{n+1}\quad
\mbox{uniformly in }X\in(\C^n)_\ast.$$

Let $a>1.$ Then we may find a neighborhood $\mathcal U$ of $\partial
D$ such that
$$a\beta_D(\zeta;X)\ge\sqrt{n+1}\kappa_D(\zeta;X),\quad\zeta\in D\cap\mathcal
U,\ X\in\C^n.$$

For $z\in D\cap\mathcal U$ and $w\in D$ ($w\neq z$) choose a smooth
curves $\gamma:[0,1]\to D\cap U$ with $\gamma(0)=z,$ $\gamma(1)=w$
and
$$a.b_D(z,w)>\inf_\gamma\int_0^1\beta_D(\gamma(t);\gamma'(t))dt.$$

Set $t_0=\sup\{t:\gamma([0,t])\subset\mathcal U\},$
$t_1=\inf\{t:\gamma([t,1])\subset\mathcal U\}$ if $w\in\mathcal U$
and $t_1=1$ otherwise. Putting $\hat z=\gamma(t_0)$ and $\hat
w=\gamma(t_1),$ it follows that
$$a.b_D(z,w)>\frac{\sqrt{n+1}}{a}[k_D(z,\hat z)+k_D(\hat w,w)]
+b_D(\hat z,\hat w)$$
$$\ge\frac{\sqrt{n+1}}{a}[k_D(z,w)-k_D(\hat
z,\hat w)].$$

If $\gamma([0,1])\subset\mathcal U,$ then
$$a.b_D(z,w)>\frac{\sqrt{n+1}}{a}k_D(z,w).$$

Otherwise, $\hat z,\hat w\not\in\mathcal U$ and then
$$ab_D(z,w)>b_D(z,\hat z)\ge\inf_{\tilde z\not\in\mathcal U} c_D(z,\tilde z)\to+\infty$$
as $z\to\partial D$ (cf. \cite[Theorem 10.2.1]{JP}).

Since $\displaystyle\sup_{\tilde z,\tilde w\not\in\mathcal
U}k_D(\tilde z,\tilde w)<+\infty,$ we get that
$$\liminf_{\substack{z\to\partial D\\z\neq
w}}\frac{k_D(z,w)}{b_D(z,w)}\ge\frac{a^2}{\sqrt{n+1}}\quad\mbox{uniformly
in }w\in D.$$

The opposite inequality
$$\liminf_{\substack{z\to\partial D\\z\neq
w}}\frac{b_D(z,w)}{k_D(z,w)}\ge\frac{\sqrt{n+1}}{a^2}\quad\mbox{(uniformly
in }w\in D)$$ follows in the same way.

To obtain \eqref{b}, it remains to let $a\to 1.$


\begin{thebibliography}{}

\bibitem{BB1} Z. M. Balogh, M. Bonk, {\it Pseudoconvexity and Gromov
hyperbolicity}, C. R. Acad. Sci. Paris S\' er. I Math. 328 (1999),
597-602.

\bibitem{BB2} Z. M. Balogh, M. Bonk, {\it Gromov hyperbolicity and
the Kobayashi metric on strictly pseudoconvex domains}, Comment.
Math. Helv. 75 (2000), 504-533.

\bibitem{For} J. E. Fornaess, {\it Strictly pseudoconvex domains in
convex domains}, Amer. J. Math. 98 (1976), 529-569.

\bibitem{Hor} L. H\"ormander, {\it $L^2$ estimates and existence theorems
for the $\overline\partial$ operator}, Acta Math. 113 (1965), 89-152.

\bibitem{JP} M. Jarnicki, P. Pflug, {\it Invariant distances and metrics
in complex analysis}, de Gruyter, 1993.

\bibitem{Lem} L. Lempert, {\it Holomorphic retracts and intrinsic metrics
in convex domains}, Analysis Mathematica 8 (1982), 257-261.

\bibitem{Ma} D. Ma, {\it Sharp estimates of the Kobayashi metric near
strongly pseudoconvex points}, Contemp. Math. 137 (1992), 329-338.

\bibitem{Nik} N. Nikolov, {\it Localization of invariant metrics},
Arch. Math. 79 (2002), 67-73.

\bibitem{NPT} N. Nikolov, P. Pflug, P. J. Thomas, {\it Upper bound
for the Lempert function of smooth domains}, Math. Z. 266 (2010),
425-430.

\bibitem{Ven} S. Venturini, {\it Comparison between the Kobayashi and
Carathe\'odory distances on strongly pseudoconvex bounded domains in
$\C^n$}, Proc. Amer. Math. Soc. 107 (1989), 725-730.

\end{thebibliography}
\end{document}